\documentclass[graybox]{svmult}

\usepackage{type1cm}          \usepackage{makeidx}          \usepackage{graphicx}         \usepackage{multicol}         \usepackage[bottom]{footmisc} 

\usepackage{newtxtext}        \usepackage[varvw]{newtxmath}       

\usepackage{cprotect}

\makeindex

\usepackage{amscd}
\usepackage{amsmath}

\newcommand{\ols}[1]{\mskip.5\thinmuskip\overline{\mskip-.5\thinmuskip {#1} \mskip-.5\thinmuskip}\mskip.5\thinmuskip} \newcommand{\olsi}[1]{\,\overline{\!{#1}}} \makeatletter
\newcommand\closure[1]{
  \tctestifnum{\count@stringtoks{#1}>1} {\ols{#1}} {\olsi{#1}} }
\long\def\count@stringtoks#1{\tc@earg\count@toks{\string#1}}
\long\def\count@toks#1{\the\numexpr-1\count@@toks#1.\tc@endcnt}
\long\def\count@@toks#1#2\tc@endcnt{+1\tc@ifempty{#2}{\relax}{\count@@toks#2\tc@endcnt}}
\def\tc@ifempty#1{\tc@testxifx{\expandafter\relax\detokenize{#1}\relax}}
\long\def\tc@earg#1#2{\expandafter#1\expandafter{#2}}
\long\def\tctestifnum#1{\tctestifcon{\ifnum#1\relax}}
\long\def\tctestifcon#1{#1\expandafter\tc@exfirst\else\expandafter\tc@exsecond\fi}
\long\def\tc@testxifx{\tc@earg\tctestifx}
\long\def\tctestifx#1{\tctestifcon{\ifx#1}}
\long\def\tc@exfirst#1#2{#1}
\long\def\tc@exsecond#1#2{#2}
\makeatother

\newcommand{\compute}[1]{\ols{#1}}
\newcommand{\extrinsic}[1]{\breve{#1}}

\usepackage{tikz}
\usetikzlibrary{patterns}

\newcommand{\Projection}{\mathscr P}
\newcommand{\DC}{\mathrm{b}}

\newcommand{\grad}{\operatorname{grad}}
\newcommand{\curl}{\operatorname{curl}}

\newcommand{\divergence}{\operatorname{div}}

\newcommand{\RT}{\bfR\bfT}

\makeatletter
\newcommand\suchthat{\@ifstar
  {\mathrel{}\middle|\mathrel{}}
  {\mid}}
\makeatother

\usepackage{accents}

\usepackage{comment}
\usepackage{url}
\usepackage[all,tips]{xy}

\usepackage{amsfonts}
\usepackage{amsmath}
\usepackage{amscd}
\usepackage{enumerate}
\usepackage{tikz}

\DeclareSymbolFont{bbold}{U}{bbold}{m}{n}
\DeclareSymbolFontAlphabet{\mathbbold}{bbold}

\newcommand{\Jacobian}{\nabla}

\newcommand{\bbR}{{\mathbb R}}

\newcommand{\bfC}{{\mathbf C}}

\newcommand{\bfF}{{\mathbf F}}

\newcommand{\bfH}{{\mathbf H}}

\newcommand{\bfL}{{\mathbf L}}

\newcommand{\bfN}{{\mathbf N}}

\newcommand{\bfR}{{\mathbf R}}

\newcommand{\bfT}{{\mathbf T}}
\newcommand{\bfU}{{\mathbf U}}
\newcommand{\bfV}{{\mathbf V}}

\newcommand{\calM}{{\mathcal M}}

\newcommand{\calP}{{\mathcal P}}

\newcommand{\calT}{{\mathcal T}}

\newcommand{\Cont}{C}
\newcommand{\Lebesgue}{L}

\newcommand{\vecCont}{\bfC}
\newcommand{\vecLebesgue}{\bfL}

\newcommand{\Manifold}{M}

\newcommand{\compManifold}{\compute{M}}

\newcommand{\Mesh}{\calT}
\newcommand{\compMesh}{\compute{\calT}}

\newcommand{\Interpolant}{I}
\newcommand{\smoothedinterpol}{Q}
\newcommand{\smoothedproj}{\pi}

\newcommand{\Ariinc}{A_h}

\newcommand{\homeo}{\Theta}

\newcommand{\metric}{g}

\newcommand{\Ned}{{\bfN\bf{e}\bf{d}}}

\begin{document}

\title*{Towards finite element exterior calculus on manifolds: commuting projections, geometric variational crimes, and approximation errors}
\author{Martin W. Licht}
\institute{École Polytechnique Fédérale de Lausanne}

\maketitle

\abstract{
 We survey recent contributions to finite element exterior calculus on manifolds and surfaces within a comprehensive formalism for the error analysis of vector-valued partial differential equations on manifolds.  Our primary focus is on uniformly bounded commuting projections on manifolds: these projections map from Sobolev de~Rham complexes onto finite element de~Rham complexes, commute with the differential operators, and satisfy uniform bounds in Lebesgue norms. They enable the Galerkin theory of Hilbert complexes for a large range of intrinsic finite element methods on manifolds.   However, these intrinsic finite element methods are generally not computable and thus primarily of theoretical interest.  This leads to our second point: estimating the geometric variational crime incurred by transitioning to computable approximate problems.  Lastly, our third point addresses how to estimate the approximation error of the intrinsic finite element method in terms of the mesh size.  If the solution is not continuous, then such an estimate is achieved via modified Cl\'ement or Scott-Zhang interpolants that facilitate a broken Bramble--Hilbert lemma. }

\section{Introduction}

We extend the mathematical theory of finite element methods for partial differential equations on surfaces and manifolds. Applications of surface finite element methods include the numerical modeling of phenomena such as biological processes on cellular membranes,  fluid dynamics within cracks, and electromagnetism on thin sheets~\cite{christiansen2002resolution,camacho2015L2,bonito2020divergence,bachini2023diffusion}. Numerical methods for partial differential equations on manifolds pertain to astrophysical applications such as modeling the star formation within relativistic regimes and black hole dynamics~\cite{arnold2000numerical,holst2001adaptive,holst2012geometric}. Whereas the numerical theory of scalar PDEs on surfaces and manifolds is well-established, numerical methods for PDEs in vector and tensor fields still demand further elaboration.  Owing to the pivotal role of differential geometry and topology in the theory of vector field PDEs, finite element exterior calculus appears to be the natural choice of formalism for such an endeavor.

Finite element exterior calculus (FEEC,~\cite{hiptmair2002finite,hiptmair2006auxiliary,arnold2010finite}) is a comprehensive mathematical framework for mixed finite element methods in vector and tensor variables through the perspective of exterior calculus. Differential complexes are central to both the theoretical and the numerical analysis of such PDEs. Geometric and topological structures are critical in understanding the stability and convergence of these methods. Exterior calculus, i.e., the calculus of differential forms, may seem challenging at first, but it is common throughout geometry and theoretical physics.  Connecting numerical analysis with exterior calculus provides access to a vast body of literature in differential geometry, functional analysis, and theoretical physics.

The following exposition summarizes recent advancements in finite element exterior calculus on manifolds. Presented within a formalism inspired by surface finite element methods, we approach numerical methods for partial differential equations on embedded surfaces as well as on Riemannian manifolds. 

First, we address the stability and quasi-optimality of mixed finite element methods for vector-valued PDEs on manifolds, which relies on commuting projections from the $\Lebesgue^{2}$ de~Rham complex onto finite element de~Rham complexes on the exact manifold.  Such commuting projections have been known for polyhedral domains~\cite{christiansen2008smoothed,licht2019smoothed} but have been developed for manifolds only recently~\cite{licht2023smoothed}.  Secondly, we address the geometric variational crime incurred by transitioning from an intrinsic finite element method, which is generally not computable and thus plays a primarily theoretical role, towards a computable finite element method.  Thirdly, we bound the quasi-optimal approximation error of the intrinsic finite element method in terms of the mesh size and the exact solution's Sobolev regularity.  This last step requires a generalized interpolant, such as of Cl\'ement-type or Scott-Zhang-type, which facilitates what has recently emerged under the term \emph{broken Bramble-Hilbert lemma}~\cite{camacho2015L2,veeser2016approximating,licht2021local,holst2023geometric}.
We address those three topics with a simple model problem of Maxwell-type without elaboration of technical details.

\section{Finite element methods for Maxwell's equations on manifolds}

Assume that $\Manifold$ is a compact oriented $3$-dimensional manifold without boundary and equipped with a Riemannian metric $\metric$. Let $\Lebesgue^{2}(\Manifold)$ and $\vecLebesgue^{2}(\Manifold)$ be the spaces of square-integrable scalar and vector fields on $\Manifold$, respectively, and     
\begin{gather*}
    H^{1}(\Manifold)
    :=
    \left\{\; U \in \Lebesgue^{2}(\Manifold) \suchthat \grad U \in \vecLebesgue^{2}(\Manifold) \;\right\}
    ,
    \\
    \bfH(\curl,\Manifold)
    :=
    \left\{\; \bfU \in \vecLebesgue^{2}(\Manifold) \suchthat \curl \bfU \in \vecLebesgue^{2}(\Manifold) \;\right\}
    .
\end{gather*}
Our physical model problem shall be the following $\curl$-$\curl$-type PDE with Lagrange multiplier~\cite{hiptmair2006auxiliary}. 
Given a vector field $\bfF \in \vecLebesgue^{2}(\Manifold)$, we search $\bfU \in \bfH(\curl,\Manifold)$ and $\zeta \in H^{1}(\Manifold)$ satisfying the following weak formulation:\footnote{For the ease of illustration, we tacitly assume that $\calM$ has first Betti number equal zero so that we can avoid the discussion of harmonic vector fields.} 
\begin{align*}
    \int_{\Manifold} \langle \curl \bfU, \curl \bfV \rangle_{\metric} + \langle \grad\zeta, \bfV \rangle_{\metric}
    &=
    \int_{\Manifold} \langle \bfF, \bfV \rangle_{\metric},
    \quad 
    &&\bfV \in \bfH(\curl,\Manifold),
    \\
    \int_{\Manifold} 
    \langle \bfU, \grad\tau \rangle_{\metric}
    - 
    \langle \zeta, \tau \rangle_{\metric} 
    &=
    0,
    \quad 
    &&\tau \in H^{1}(\Manifold)
    .
\end{align*}
Let $\Mesh$ be a triangulation of the manifold $\Manifold$ such that each cell $T \in \Mesh$ is diffeomorphic to the reference tetrahedron $\hat T$ along a reference diffeomorphism $\imath_T : \hat T \to T$. 
We emphasize that $\imath_{T}$ is part of the definition of each cell, and require these transformations to be compatible along interfaces. 
This intrinsic triangulation $\Mesh$ is first and foremost of theoretical relevance; its relation with computational triangulations will be discussed later. For now, we mention that in practice we have uniformly bounded constants $C_{m} > 0$ such that for each $T \in \Mesh$ with $\metric$-diameter $h_T$: 
\begin{align*}
    \| \Jacobian^{m} \imath_T \|_{\Lebesgue^{\infty}(\hat T)}
    \leq 
    C_{m} h_T^{m}
    ,
    \quad 
    \| \Jacobian^{m} \imath_T^{-1} \|_{\Lebesgue^{\infty}(T)}
    \leq 
    C_{m} h_T^{-m}
    .
\end{align*}
With that preparation in place, 
we define finite element spaces intrinsically by piecewise pullback from the reference tetrahedron and then imposing the relevant continuity conditions along cell interfaces. 
For instance, the intrinsic N\'ed\'elec space $\Ned_{r}(\Mesh)$ is the piecewise transform of the N\'ed\'elec vector fields of polynomial degree $r \geq 0$ with tangential continuity imposed along faces.
This leads to an intrinsic finite element formulation: 
we search $( \bfU_{h}, \zeta_{h} ) \in \Ned_{r}(\Mesh) \times \calP_{r+1}(\Mesh)$ satisfying:
\begin{align*}
    \int_{\Manifold} \langle \curl \bfU_{h}, \curl \bfV_{h} \rangle_{\metric} + \langle \grad\zeta_{h}, \bfV_{h} \rangle_{\metric}
    &=
    \int_{\Manifold} \langle \bfF, \bfV_{h} \rangle_{\metric},
    \quad 
    &&\bfV_{h} \in \Ned_{r}(\Mesh),
    \\
    \int_{\Manifold} 
    \langle \bfU_{h}, \grad\tau_{h} \rangle_{\metric}
    - 
    \langle \zeta_{h}, \tau_{h} \rangle_{\metric} 
    &=
    0,
    \quad 
    &&\tau_{h} \in \calP_{r+1}(\Mesh)
    .
\end{align*}
The stability and convergence of this finite element method follows easily within the Galerkin theory of Hilbert complexes~\cite{arnold2010finite} once we have a commuting projection.

When discussing partial differential equations in vector fields, it is advisable to first recall the Sobolev de~Rham complex:
\begin{gather*}
    \begin{CD}
        H^{1}(\Manifold{})
        @>{\grad}>>
        \bfH(\curl,\Manifold{})
        @>{\curl}>>
        \bfH(\divergence,\Manifold{})
        @>{\divergence}>>
        \Lebesgue^{2}(\Manifold{})
        .
    \end{CD}
\end{gather*}
Analogously, the discussion of finite element methods for PDEs in vector fields is best understood on the background of finite element de~Rham complexes. 
We may consider, for instance, the following finite element de~Rham complex:
\begin{gather*}
    \begin{CD}
        \calP_{r+1}(\Mesh{})
        @>{\grad}>>
        \Ned_{r}(\Mesh{})
        @>{\curl}>>
        \RT_{r}(\Mesh{})
        @>{\divergence}>>
        \calP_{r,\DC}(\Mesh{})
        .
    \end{CD}
\end{gather*}
Here, the spaces are the Lagrange space $\calP_{r+1}(\Mesh{})$, the N\'ed\'elec space $\Ned_{r}(\Mesh{})$ and Raviart-Thomas space $\RT_{r}(\Mesh{})$, and the broken Lagrange space $\calP_{r,\DC}(\Mesh{})$.
In order to utilize the full scope of the Galerkin theory of Hilbert complexes, 
we want $L^{2}$-bounded idempotent mappings that make the following diagram commute:
\begin{gather*}
    \begin{CD}
        H^{1}(\Manifold{})
        @>{\grad}>>
        \bfH(\curl,\Manifold{})
        @>{\curl}>>
        \bfH(\divergence,\Manifold{})
        @>{\divergence}>>
        \Lebesgue^{2}(\Manifold{})
        \\
        @V{\smoothedproj_{h}^{0}}VV
        @V{\smoothedproj_{h}^{1}}VV
        @V{\smoothedproj_{h}^{2}}VV
        @V{\smoothedproj_{h}^{3}}VV
        \\
        \calP_{r+1}(\Mesh{})
        @>{\grad}>>
        \Ned_{r}(\Mesh{})
        @>{\curl}>>
        \RT_{r}(\Mesh{})
        @>{\divergence}>>
        \calP_{r,\DC}(\Mesh{})
        . 
    \end{CD}
\end{gather*}
Given such projections, the Galerkin theory of Hilbert complexes establishes stability and quasi-optimality of the intrinsic finite element method:
$\| \bfU_{h} \|_{\bfH(\curl,\Manifold{})} + \| \zeta_{h} \|_{H^{1}(\Manifold{})} \leq C \| \bfF \|_{\vecLebesgue^{2}(\Manifold{})}$ and  
\begin{align*}
&
 \| \bfU - \bfU_{h} \|_{\bfH(\curl,\Manifold{})} 
 +
 \| \zeta - \zeta_{h} \|_{H^{1}(\Manifold{})} 
 \\
 &\qquad\qquad 
 \leq 
 C
 \inf_{ (\bfV_{h},\tau_{h}) \in \Ned_{r}(\Mesh{}) \times \calP_{r+1}(\Mesh{}) }
 \| \bfU - \bfV_{h} \|_{\bfH(\curl,\Manifold{})} 
 +
 \| \zeta - \tau_{h} \|_{H^{1}(\Manifold{})} 
 .
\end{align*}
The commuting projection is composed of several operators. We refer to~\cite{licht2023smoothed} for the full discussion. 
First, a smoothing operator maps from the $\Lebesgue^{2}$ de~Rham complex onto the smooth de~Rham complex, which leads to a commuting diagram:
\begin{gather*}
    \begin{CD}
        H^{1}(\Manifold{})
        @>{\grad}>>
        \bfH(\curl,\Manifold{})
        @>{\curl}>>
        \bfH(\divergence,\Manifold{})
        @>{\divergence}>>
        \Lebesgue^{2}(\Manifold{})
        \\
        @V{R^{0}}VV
        @V{R^{1}}VV
        @V{R^{2}}VV
        @V{R^{3}}VV
        \\
        \Cont^{\infty}(\Manifold)
        @>{\grad}>>
        \vecCont^{\infty}(\Manifold)
        @>{\curl}>>
        \vecCont^{\infty}(\Manifold)
        @>{\divergence}>>
        \Cont^{\infty}(\Manifold)
        . 
    \end{CD}
\end{gather*}
Here, we expand upon an old idea by de~Rham~\cite{derham1984differentiable}, defining a localized smoothing operator within each coordinate chart. Away from the boundary of the coordinate chart, the local operator appears like convolution with a smooth mollifier, but the convolution radius shrinks to zero towards the boundary. The global smoothing operator is then the composition of these local operators: a minimum amount of smoothing is applied everywhere due to some local smoothing operator. 

The second stage of the smoothed projection is the canonical finite interpolant. 
These interpolants fill up another commuting diagram:
\begin{gather*}
    \begin{CD}
        \Cont^{\infty}(\Manifold)
        @>{\grad}>>
        \vecCont^{\infty}(\Manifold)
        @>{\curl}>>
        \vecCont^{\infty}(\Manifold)
        @>{\divergence}>>
        \Cont^{\infty}(\Manifold)
        \\
        @V{\Interpolant_{h}^{0}}VV
        @V{\Interpolant_{h}^{1}}VV
        @V{\Interpolant_{h}^{2}}VV
        @V{\Interpolant_{h}^{3}}VV
        \\
        \calP_{r+1}(\Mesh)
        @>{\grad}>>
        \Ned_{r}(\Mesh)
        @>{\curl}>>
        \RT_{r}(\Mesh)
        @>{\divergence}>>
        \calP_{r,\DC}(\Mesh)
        . 
    \end{CD}
\end{gather*}
The commuting quasi-interpolant $\smoothedinterpol_{h}^{k} = \Interpolant_{h}^{k} R^{k}$ is the composition of the canonical interpolant with the smoothing operators. While they commute with the differential operators and satisfy uniform bounds in Lebesgue norms, they are generally not idempotent. However, we can uniformly bound their difference from the identity mapping over the finite element spaces. As a consequence, there exists a (generally non-local) operator that reconstructs the original vector field and which is part of another commuting diagram:
\begin{gather*}
    \begin{CD}
        \calP_{r+1}(\Mesh)
        @>{\grad}>>
        \Ned_{r}(\Mesh)
        @>{\curl}>>
        \RT_{r}(\Mesh)
        @>{\divergence}>>
        \calP_{r,\DC}(\Mesh)
        \\
        @V{{J}_{h}^{0}}VV
        @V{{J}_{h}^{1}}VV
        @V{{J}_{h}^{2}}VV
        @V{{J}_{h}^{3}}VV
        \\
        \calP_{r+1}(\Mesh)
        @>{\grad}>>
        \Ned_{r}(\Mesh)
        @>{\curl}>>
        \RT_{r}(\Mesh)
        @>{\divergence}>>
        \calP_{r,\DC}(\Mesh)
        . 
    \end{CD}
\end{gather*}
The composition of these stages defines the smoothed projection $\smoothedproj_{h}^{k} = J_{h}^{k} \Interpolant_{h}^{k} R^{k}$.
The literature on mixed finite element methods knows different commuting projections for different geometric settings, and with different continuity and locality properties~\cite{falk2014local,ern2016mollification,licht2019smoothed,licht2019mixed,arnold2021local,ern2022equivalence}, 
but the extension to manifolds~\cite{licht2023smoothed} is a recent contribution.

\section{Computational Problem and Geometric Error} 

Our model problem and the intrinsic finite element method are posed on the exact physical manifold $\Manifold$ with Riemannian metric $\metric$. However, this exact geometric ambient is not accessible in practical computations. We assess the effect of the \emph{geometric variational crime} incurred by transitioning to an approximate geometric ambient. For that purpose, we introduce the approximate computational manifold $\compManifold$, which may generally be different from $\Manifold$,
and which is equipped with a computational Riemannian metric $\compute{\metric}$.
The manifold $\compManifold$ carries a triangulation $\compMesh$.

The relation between $\Manifold$ and $\compManifold$ is described by a homeomorphism $\homeo : \compManifold \to \Manifold$ with the following properties:
(i) $\homeo$ is a bi-Lipschitz mapping, with respect to the metric space induced by the Riemannian tensors,
and (ii) given any cell $T \in \compMesh$, the restriction $\homeo_{|T}$ is a diffeomorphism.
We emphasize that $(\compManifold,\compute{\metric})$ together with $\compMesh$ are computable, 
whereas $\homeo$ and the pair $(\Manifold,\metric)$ are generally not. 
We also remark that the intrinsic triangulation $\Mesh$ can simply be defined as the image of the computable triangulation $\compMesh$ along the homeomorphism $\homeo$. 

\begin{example}
 The homeomorphism $\homeo : \compManifold \to \Manifold$ is explicitly computable
 if $\Manifold$ is explicitly described via coordinate charts, 
 as is the underlying assumption of some astrophysics codes~\cite{lindblom2013solving}. 
 Then $\compute{\metric}$ is taken as the interpolation of $\extrinsic{\metric}$ on $\compMesh$. 
\end{example}
\begin{example}
 In the setting of surface finite element methods~\cite{camacho2015L2,bonito2020divergence},
 the manifold $\Manifold$ is a hypersurface in $\bbR^{n}$ and $\compManifold$ is a piecewise polynomial approximate surface, defined via some interpolation procedure. 
 The Euclidean ambient metric induces the exact metric $\metric$ and the computational metric $\compute{\metric}$.
 The closest-point projection onto $\Manifold$ defines a homeomorphism $\homeo : \compManifold \to \Manifold$,
 which is generally not computable. 
\end{example}

The following finite element method is implementable in practice: 
given a vector field $\compute{\bfF} \in \vecLebesgue^{2}(\compManifold)$, 
we search $(\compute\bfU_{h},\compute\zeta_{h}) \in \Ned_{r}(\compMesh) \times \calP_{r+1}(\compMesh)$ satisfying 
    \begin{align*}
    \int_{\compManifold} \langle \curl \compute\bfU_{h}, \curl \compute\bfV_{h} \rangle_{\compute\metric} 
    + 
    \langle \grad\compute\zeta_{h}, \compute\bfV_{h} \rangle_{\compute\metric}
    &=
    \int_{\compManifold} \langle \compute\bfF, \compute\bfV_{h} \rangle_{\compute\metric}, \quad 
    &&\compute\bfV_{h} \in \Ned_{r}(\compMesh),
    \\
    \int_{\compManifold} 
    \langle \compute\bfU_{h}, \grad\compute\tau_{h} \rangle_{\compute\metric}
    - 
    \langle \compute\zeta_{h}, \compute\tau_{h} \rangle_{\compute\metric} 
    &=
    0,
    \quad 
    &&\compute\tau_{h} \in \calP_{r+1}(\compMesh)
    .
\end{align*}
Comparing the computational finite element solution $\compute\bfU_{h}$ with the exact solution $\bfU$
is facilitated via an intermediate problem. 
We let $\extrinsic{\metric}$ be the pullback of the metric on $\Manifold$ onto the computational manifold $\compManifold$. 
By construction, $\homeo : (\compManifold,\extrinsic{\metric}) \to (\Manifold,\metric)$ is an isometry of Riemannian manifolds.
The intermediate problem of finding $(\extrinsic\bfU_{h},\extrinsic\zeta_{h}) \in \Ned_{r}(\compMesh) \times \calP_{r+1}(\compMesh)$ satisfying 
\begin{align*}
    \int_{\compManifold} 
    \langle \curl \extrinsic\bfU_{h}, \curl \compute\bfV_{h} \rangle_{\extrinsic\metric} 
    + 
    \langle \grad\extrinsic\zeta_{h}, \compute\bfV_{h} \rangle_{\extrinsic\metric}
    &=
    \int_{\compManifold} \langle \homeo^{\ast}\bfF, \compute\bfV_{h} \rangle_{\extrinsic\metric},
    \quad 
    &&\compute\bfV_{h} \in \Ned_{r}(\compMesh),
    \\
    \int_{\compManifold} 
    \langle \extrinsic\bfU_{h}, \grad\compute\tau_{h} \rangle_{\extrinsic\metric}
    - 
    \langle \extrinsic\zeta_{h}, \compute\tau_{h} \rangle_{\extrinsic\metric} 
    &=
    0,
    \quad 
    &&\compute\tau_{h} \in \calP_{r+1}(\compMesh)
\end{align*}
is equivalent to the intrinsic finite element problem, i.e., $\extrinsic\bfU_{h} = \homeo^{\ast} \bfU_{h}$ and $\extrinsic\zeta_{h} = \homeo^{\ast} \zeta_{h}$.
We focus on the vector variable and estimate: 
\begin{align*}
 \| \homeo^{\ast}\bfU - \compute\bfU_{h} \|_{\bfH(\curl,\compManifold,\compute\metric)}
 &\leq 
 \| \homeo^{\ast}\bfU - \extrinsic\bfU_{h} \|_{\bfH(\curl,\compManifold,\compute\metric)}
 +
 \| \extrinsic\bfU_{h} - \compute\bfU_{h} \|_{\bfH(\curl,\compManifold,\compute\metric)}
.
\end{align*}
The identity mapping gives rise to a linear isomorphism $\Ariinc$ of finite element differential complexes that only differ by their respective Hilbert space structures:
\begin{gather*}
    \begin{CD}
        \calP_{r+1}(\compMesh,\extrinsic{\metric})
        @>{\grad}>>
        \Ned_{r}(\compMesh,\extrinsic{\metric})
        @>{\curl}>>
        \RT_{r}(\compMesh,\extrinsic{\metric})
        @>{\divergence}>>
        \calP_{r,\DC}(\compMesh,\extrinsic{\metric})
        \\
        @A{\Ariinc}AA @A{\Ariinc}AA @A{\Ariinc}AA @A{\Ariinc}AA
        \\
        \calP_{r+1}(\compMesh,\compute{\metric})
        @>{\grad}>>
        \Ned_{r}(\compMesh,\compute{\metric})
        @>{\curl}>>
        \RT_{r}(\compMesh,\compute{\metric})
        @>{\divergence}>>
        \calP_{r,\DC}(\compMesh,\compute{\metric})
        . 
    \end{CD}
\end{gather*}
We now express the effect of the variational crime via a variant of Theorem~3.10 of~\cite{holst2012geometric},
which leads to 
\begin{align*}
    &
    \| \extrinsic{{\bfU}}_h - \Ariinc \compute{{\bfU}}_h \|_{\bfH(\curl,\compManifold,\compute{\metric})}
    \leq
    C 
    \left(
      \| \compute{\bfF} - \Ariinc^{\ast} \theta^{\ast}{\bfF} \|_{\vecLebesgue^{2}(\compManifold,\compute{\metric})} 
      + 
      \| I - \Ariinc^\ast\Ariinc \| 
      \| \bfF \|_{\vecLebesgue^{2}(\compManifold,\compute{\metric})}
    \right).
\end{align*}
The difference $\compute{\bfF} - \Ariinc^{\ast} \theta^{\ast}{\bfF}$ depends on implementation details and we refer to Theorem~3.12 of~\cite{holst2012geometric} for some theoretical results. 
The analysis of the geometric error depends on bounding the quantity 
\begin{align*}
    \left\| I - \Ariinc^\ast\Ariinc \right\|
    &=
    \sup_{ \compute{\bfU} \in \vecLebesgue^{2}(\compManifold) }
    \left| 
        1 - 
        \textstyle{
        \int_{\compManifold} \langle \compute{\bfU}, \compute{\bfU} \rangle_{\compute{\metric}} 
        \big/
        \int_{\compManifold} \langle \compute{\bfU}, \compute{\bfU} \rangle_{\extrinsic{\metric}} 
        }
\right|
    .
\end{align*}
Estimating this depends, again, on implementation details.
Intuitively, this measures in how far the computational metric $\compute{\metric}$ differs from the pullback metric $\extrinsic{\metric}$, or in how far $\homeo : (\compManifold,\compute{\metric}) \to (\Manifold,\metric)$ differs from being an isometry.

\section{Approximation error}

Lastly, we explain recent results in the approximation theory of finite element vector analysis. 
Any complete error analysis of the finite element method requires estimating $\| \extrinsic\bfU_{h} - \homeo^{\ast}\bfU \|_{\bfH(\curl,\compManifold,\compute\metric)}$. Here, we have transferred the physical solution onto the computational manifold. 
Alternatively, we transfer the computational finite element solution onto the physical manifold and estimate 
$\homeo^{-\ast}\extrinsic\bfU_{h} - \bfU$. 
In either case, we want to derive error estimates in terms of mesh size and the solution regularity.

Let $\bfH^{m}(\Manifold)$ be the $m$-th order Sobolev space of vector fields on $\Manifold$, that is, the space of vector fields whose covariant derivatives up to order $m$ are square-integrable. 
This is a Banach space with norm 
$\| \bfU \|_{\bfH^{m}(\Manifold)} := \sum_{k=0}^{m} \| \nabla^{k} \bfU \|_{\Lebesgue^{2}(\Manifold)}$.
Even if $\bfU \in \bfH^{m}(\Manifold)$, we cannot directly exploit this regularity for a Bramble--Hilbert-type estimate because there is no intrinsic notion of polynomials on manifolds.
Instead, we study the pullback $\homeo^{\ast}\bfU$ on the computational manifold. 
However, $\homeo$ is only globally Lipschitz, and the computational manifold $\compManifold$ has only a piecewise smooth structure. 
As a consequence, we generally only have $\homeo^{\ast}\bfU \in \bfH(\curl,\compManifold)$ and $\homeo^{\ast}\bfU_{|T} \in \bfH^{m}(T)$ for each $T \in \compMesh$. 
The situation is easy if $\bfU$ is continuous: 
then $\homeo^{\ast}\bfU$ has the necessary interelement continuity so that the canonical N\'ed\'elec interpolant $I_{\Ned}$ is applicable. Together with the standard Bramble--Hilbert lemma, we estimate 
\begin{align*}
 \| \homeo^{\ast}\bfU - \extrinsic\bfU_{h} \|_{\bfH(\curl,\compManifold,\compute\metric)}
 &
 \leq 
 C
 \| \homeo^{\ast}\bfU - I_{\Ned} \homeo^{\ast}\bfU \|_{\bfH(\curl,\compManifold,\compute\metric)}
\leq 
 C
 \sum_{T \in \Mesh} h_T^{m} \| \bfU \|_{\bfH^{m}(T)}.
\end{align*}
While that may suffice for our simple model problem, 
solutions of more complicated PDEs and with discontinuous coefficients may have insufficient regularity for the canonical interpolant. 
However, a generalized Scott-Zhang projection~\cite{licht2021local} or Cl\'ement-Ern-Guermond projection for manifolds 
\begin{gather*}
    \Projection_{\Ned} : \bfH(\curl,\compManifold) \rightarrow \Ned_{r}(\compManifold),
\end{gather*}
enables the following error estimate:
for any $m \in [0,r+1]$, $l \in [0,r+1]$, all tetrahedra $T \in \compManifold$, and all $\compute{\bfV} \in \bfH(\curl,\compManifold)$, it holds that 
\begin{gather*}
    \| \compute{\bfV} - \Projection_{\Ned} \compute{\bfV} \|_{\vecLebesgue^{2}(T)}
    \leq 
    C        
    \sum_{ \substack{ T' \in \compManifold \\ \dim(T') = 3 \\ T' \cap T \neq \emptyset } }
    h_T^{m} \| \compute{\bfV} \|_{\bfH^{m}(T')}
    +
    h_T^{l+1} \| \curl \compute{\bfV} \|_{\bfH^{l}(T')}
    .
\end{gather*}
Here, $C > 0$ depends only on $r$ and the mesh regularity. This last inequality is also known as \emph{broken Bramble-Hilbert lemma} as it only relies on the \emph{broken} or \emph{piecewise} regularity of $\homeo^{\ast}\bfU$. In particular, the regularity required is less than that of the Cl\'ement interpolant, whose error bounds would still require a modicum of Sobolev regularity even across cells.  By contrast, the projection $\Projection_{\Ned}$ enables error estimates even in situations of lowest regularity.


\begin{thebibliography}{10}

\bibitem{arnold2010finite}
{\sc D.~Arnold, R.~Falk, and R.~Winther}, {\em Finite element exterior
  calculus: from {H}odge theory to numerical stability}, Bulletin of the
  American mathematical society, 47 (2010), pp.~281--354.

\bibitem{arnold2021local}
{\sc D.~Arnold and J.~Guzm{\'a}n}, {\em Local ${L}^2$-bounded commuting
  projections in {FEEC}}, ESAIM: Mathematical Modelling and Numerical Analysis,
  55 (2021), pp.~2169--2184.

\bibitem{arnold2000numerical}
{\sc D.~N. Arnold}, {\em Numerical problems in general relativity}, Numerical
  Mathematics and Advanced Applications (P. Neittaanmki, T. Tiihonen, and P.
  Tarvainen, eds.), World Scientific,  (2000), pp.~3--15.

\bibitem{bachini2023diffusion}
{\sc E.~Bachini, P.~Brandner, T.~Jankuhn, M.~Nestler, S.~Praetorius,
  A.~Reusken, and A.~Voigt}, {\em Diffusion of tangential tensor fields:
  numerical issues and influence of geometric properties}, Journal of Numerical
  Mathematics,  (2023).

\bibitem{bonito2020divergence}
{\sc A.~Bonito, A.~Demlow, and M.~Licht}, {\em A divergence-conforming finite
  element method for the surface {S}tokes equation}, SIAM Journal on Numerical
  Analysis, 58 (2020), pp.~2764--2798.

\bibitem{camacho2015L2}
{\sc F.~Camacho and A.~Demlow}, {\em ${L}^2$ and pointwise a posteriori error
  estimates for {FEM} for elliptic {PDE}s on surfaces}, IMA Journal of
  Numerical Analysis, 35 (2015), pp.~1199--1227.

\bibitem{christiansen2008smoothed}
{\sc S.~Christiansen and R.~Winther}, {\em Smoothed projections in finite
  element exterior calculus}, Mathematics of Computation, 77 (2008),
  pp.~813--829.

\bibitem{christiansen2002resolution}
{\sc S.~H. Christiansen}, {\em R{\'e}solution des {\'e}quations int{\'e}grales
  pour la diffraction d'ondes acoustiques et
  {\'e}lectromagn{\'e}tiques-Stabilisation d'algorithmes it{\'e}ratifs et
  aspects de l'analyse num{\'e}rique}, PhD thesis, Ecole Polytechnique X, 2002.

\bibitem{derham1984differentiable}
{\sc G.~De~Rham}, {\em Differentiable manifolds: forms, currents, harmonic
  forms}, vol.~266, Springer Science \& Business Media, 2012.

\bibitem{ern2022equivalence}
{\sc A.~Ern, T.~Gudi, I.~Smears, and M.~Vohral{\'\i}k}, {\em Equivalence of
  local-and global-best approximations, a simple stable local commuting
  projector, and optimal $hp$ approximation estimates in ${H}(div)$}, IMA
  Journal of Numerical Analysis, 42 (2022), pp.~1023--1049.

\bibitem{ern2016mollification}
{\sc A.~Ern and J.-L. Guermond}, {\em Mollification in strongly {L}ipschitz
  domains with application to continuous and discrete de {R}ham complexes},
  Computational Methods in Applied Mathematics, 16 (2016), pp.~51--75.

\bibitem{falk2014local}
{\sc R.~Falk and R.~Winther}, {\em Local bounded cochain projections},
  Mathematics of Computation, 83 (2014), pp.~2631--2656.

\bibitem{licht2021local}
{\sc E.~Gawlik, M.~J. Holst, and M.~W. Licht}, {\em Local finite element
  approximation of {S}obolev differential forms}, ESAIM: M2NA, 55 (2021),
  pp.~2075--2099.

\bibitem{hiptmair2002finite}
{\sc R.~Hiptmair}, {\em Finite elements in computational electromagnetism},
  Acta Numerica, 11 (2002), pp.~237--339.

\bibitem{hiptmair2006auxiliary}
{\sc R.~Hiptmair, G.~AU, Widmer, and J.~Zou}, {\em Auxiliary space
  preconditioning in ${H}_0(curl; \omega)$}, Numerische Mathematik, 103 (2006),
  pp.~435--459.

\bibitem{holst2001adaptive}
{\sc M.~Holst}, {\em Adaptive numerical treatment of elliptic systems on
  manifolds}, Advances in Computational Mathematics, 15 (2001), pp.~139--191.

\bibitem{holst2023geometric}
{\sc M.~Holst and M.~Licht}, {\em Geometric transformation of finite element
  methods: Theory and applications}, Applied Numerical Mathematics, 192 (2023),
  pp.~389--413.

\bibitem{holst2012geometric}
{\sc M.~Holst and A.~Stern}, {\em Geometric variational crimes: {H}ilbert
  complexes, finite element exterior calculus, and problems on hypersurfaces},
  Foundations of Computational Mathematics, 12 (2012), pp.~263--293.

\bibitem{licht2023smoothed}
{\sc M.~Licht}, {\em Smoothed projections in finite element exterior calculus
  over manifolds}.
\newblock Arxiv preprint: https://arxiv.org/abs/2310.14276v1.

\bibitem{licht2019mixed}
\leavevmode\vrule height 2pt depth -1.6pt width 23pt, {\em Smoothed projections
  and mixed boundary conditions}, Mathematics of Computation, 88 (2019),
  pp.~607--635.

\bibitem{licht2019smoothed}
\leavevmode\vrule height 2pt depth -1.6pt width 23pt, {\em Smoothed projections
  over weakly {L}ipschitz domains}, Mathematics of Computation, 88 (2019),
  pp.~179--210.

\bibitem{lindblom2013solving}
{\sc L.~Lindblom and B.~Szil{\'a}gyi}, {\em Solving partial differential
  equations numerically on manifolds with arbitrary spatial topologies},
  Journal of Computational Physics, 243 (2013), pp.~151--175.

\bibitem{veeser2016approximating}
{\sc A.~Veeser}, {\em Approximating gradients with continuous piecewise
  polynomial functions}, Foundations of Computational Mathematics, 16 (2016),
  pp.~723--750.

\end{thebibliography}
\end{document}